\documentclass[a4paper,twoside,12pt]{article}
\usepackage{amsmath,amssymb,amsthm}
\usepackage{bbm}
\usepackage{bm}
\usepackage{verbatim}
\usepackage{stmaryrd}
\usepackage{graphicx,color}
\usepackage{mathtools}
\usepackage{titlesec}
\usepackage{titletoc}
\usepackage{xparse}
\usepackage{hyperref}
\usepackage{slashed}
\usepackage{enumitem}

\usepackage{indentfirst}

\titleformat{name=\section,numberless}{\large\bfseries}{}{0pt}{}
\titleformat*{\section}{\large\bfseries}
\titleformat*{\subsection}{\normalsize\bfseries}

\DeclareMathOperator{\RE}{Re}
\DeclareMathOperator{\IM}{Im}
\DeclareMathOperator{\image}{image}
\DeclareMathOperator{\SPAN}{span}
\DeclareMathOperator{\dist}{dist}

\newcommand{\eps}{\varepsilon}

\newcommand{\dd}{\mathrm{d}}
\newcommand{\p}{\partial}

\newcommand{\tn}[1]{{\vert\kern-0.25ex\vert\kern-0.25ex\vert #1 \vert\kern-0.25ex\vert\kern-0.25ex\vert}}

\renewcommand{\epsilon}{USE eps INSTEAD}

\newcommand{\C}{\mathbbm{C}}
\newcommand{\R}{\mathbbm{R}}

\newcommand{\ZXXX}{\mathbbm{Z}}
\newcommand{\allint}{\ZXXX}
\newcommand{\nonnegint}{\mathbbm{N}_0}

\titlecontents{section}[4.2em]{\addvspace{0pt}}
{\contentslabel{1.5em}}
{}{\titlerule*[0.3pc]{.}\contentspage\rule{3.5em}{0pt}}

\renewcommand{\subset}{\subseteq}

\newcommand{\lemmaref}[1]{Lemma \ref{lemma:#1}}
\newcommand{\theoremref}[1]{Theorem \ref{theorem:#1}}

\newcommand{\ip}[2]{\langle #1,#2\rangle}
\newcommand{\ltwo}[1]{\|#1\|}

\newcommand{\op}{\mathbf{D}}
\newcommand{\dom}{\Omega}
\newcommand{\uc}{\smash{\widehat{\dom}}}
\newcommand{\x}[1]{x^{#1}}

\newcommand{\na}[1]{|A|_{#1}}
\newcommand{\nb}[1]{|B+zA^0|_{#1}}

\newcommand{\rz}{z_{\ast}}
\newcommand{\cc}{R}

\newcommand{\rr}[1]{q_{#1}}
\newcommand{\rrrz}{\rr{1\ast}}%

\newcommand{\ret}{\widehat{\op}^{-1}_{\textup{ret}}}

\newcommand{\banach}{V}

\newcommand{\auxi}{\Xi}

\newcommand{\semi}[1]{\|#1\|_{\textup{nc}}}

\newcommand{\spek}{\sigma}
\newcommand{\pspek}{\sigma_{\textup{eigenvalue}}}

\newtheoremstyle{mytheoremstyle}
{14pt}
{14pt}
{\itshape}
{20pt}
{\bfseries}
{.}
{.5em}
{}

\newtheoremstyle{myremarkstyle}
{12pt}
{12pt}
{}
{20pt}
{\itshape}
{.}
{.5em}
{}

\theoremstyle{mytheoremstyle}

\newtheorem{theorem}{Theorem}[section]
\newtheorem{lemma}[theorem]{Lemma}
\newtheorem{corollary}[theorem]{Corollary}

\theoremstyle{myremarkstyle}

\newtheorem{remark}[theorem]{Remark}
\newtheorem{example}[theorem]{Example}
\newtheorem{counterexample}[theorem]{Counterexample}
\newtheorem*{PROOF}{Proof}

\newcommand{\xx}{w}

\newcommand{\step}{\vskip 4mm}

\newcommand{\someletter}{\mathcal{X}}
\newcommand{\fst}{\someletter'}
\newcommand{\snd}{\someletter''}
\newcommand{\trd}{\someletter}
\newcommand{\smoo}{\mathcal{C}^{\infty}}
\newcommand{\smooc}{\widehat{\mathcal{C}}^{\infty}}

\newcommand{\infdiff}{$\infty$-differentiable}

\newcommand{\ban}[1]{\mathcal{H}_{#1}}
\newcommand{\banc}[2]{\widehat{\mathcal{H}}_{1,#1,#2}}
\newcommand{\fr}{\mathbf{F}}

\begin{document}

\noindent{\bf\Large
Finite codimension stability of some\\
\rule{0pt}{18pt}time-periodic hyperbolic equations\\
\rule{0pt}{18pt}(via compact resolvents)\par}

\noindent\rule{0pt}{30pt}{\bf
Michael Reiterer}
\vskip 10mm
\noindent {\bf Abstract:}
We identify a class of
time-periodic linear symmetric hyperbolic equations
that are finite codimension stable,
because an associated operator
has compact resolvent, sufficiently
far to the right in the complex plane.
This paper is an attempt to
capture abstractly
the observation in numerical general relativity
that some discretely self-similar spacetimes,
such as Choptuik's critical spacetime,
are finite codimension stable.
\vskip 8mm

\setcounter{tocdepth}{1}
\tableofcontents


\newcommand{\selfsim}{\mathbf{T}}
\section{Introduction}\label{sec:intro}

\subsection*{Motivation}

This paper is motivated by 
a question in general relativity: \emph{%
How does one explain the numerical observation
that some discretely self-similar spacetimes,
such as Choptuik's,
are stable with small finite codimension \cite{choptuik1993}, \cite{GM}?}
One might hope for an explanation that is abstract,
and neither requires general relativity
nor depends on details of spacetimes.
\step
It is plausible that many aspects of stability
are already present at the linearized level,
and this paper is about linear equations only.
\step
We identify simple abstract assumptions
for a time-periodic linear symmetric hyperbolic equation that imply
that it is finite codimension stable,
roughly: if in the space of solutions one mods out the decaying solutions,
only finitely many independent
growing solutions are left.
This is due to an effect at high frequencies,
whereas to actually determine the codimension,
one would need a more detailed understanding of low frequencies\footnote{%
Our abstract assumptions are consistent with all
finite values of the codimension.}\textsuperscript{,}\footnote{%
Better understanding low frequencies will often be an
effectively finite-dimensional problem, amenable to (say)
rigorous computer-assisted study.
Whether there is a more conceptual, abstract
approach to low frequencies is an interesting question.
}.
\step

\newcommand{\intpic}{scaling-pic}
\newcommand{\wrkpic}{translation-pic}
This introduction is informal;
the logical development starts in Section \ref{sec:aa}.
We are on an infinite solid cylinder,
$\R \times \text{(closed $n$-dim ball)}$.
It is useful to have two different but
diffeomorphic pictures of this cylinder:
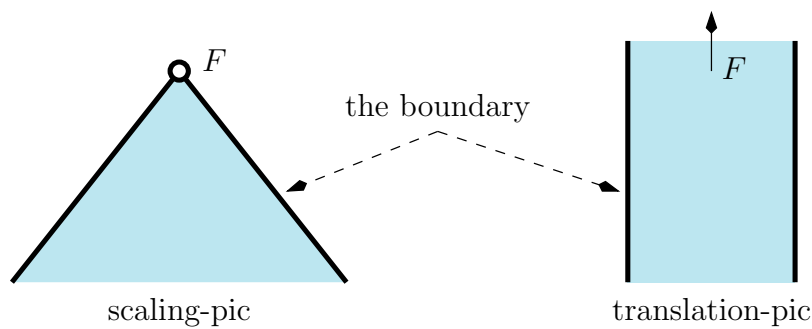
\begin{center}
\input{two.pstex_t}
\end{center}
In both, time increases upwards;
space is horizontal; horizontal cross-sections are closed $n$-dim balls.
Here $F \notin \text{cylinder}$ is the future limit point.
Let $\selfsim$ be a fixed self-diffeomorphism of the cylinder,
equivalently given:
\begin{itemize}
\item In the \intpic, by a scaling about and towards $F$.
\item In the \wrkpic, by a translation upwards.
\end{itemize}

\newcommand{\ncomp}{\text{cylinder}\to\C^N}
Consider a linear symmetric hyperbolic operator
$\widehat{\op}$ on the cylinder.
That is, an $N\times N$ matrix whose entries are
first order
differential operators,
subject to a standard condition on its principal part.
We assume:
\begin{itemize}
\item \emph{$\selfsim$-periodicity:}
$\widehat{\op}$ commutes with composition by $\selfsim$,
that is,
for all
functions $\widehat{u} : \ncomp$ one has
$(\widehat{\op}\widehat{u})\circ \selfsim
= \widehat{\op}(\widehat{u}\circ \selfsim)$.
\item \emph{Causal independence:}
$\widehat{\op}$ allows
information to flow out through the boundary,
but no information to flow in,
as time increases.
\end{itemize}
This operator $\widehat{\op}$ could arise in general relativity
from linearizing the Einstein equations about a
spacetime for which $\selfsim$ is a constant rescaling of the metric\footnote{%
Beware that the linearized Einstein equations
are not symmetric hyperbolic out of the box.
Rather, the space of solutions to say the
linearized vacuum Einstein equations,
about a background solution,
is the first cohomology of a differential graded Lie algebra
\cite{gradedliealgebra}.
The equations can be made symmetric hyperbolic through gauge-fixing.
},
aka a discretely self-similar spacetime,
or it could arise somewhere else.
\step
For every nice function $\widehat{f}: \ncomp$
with compact support in time,
let $\ret \widehat{f} : \ncomp$ be the unique solution to
\[
\widehat{\op}(\ret \widehat{f}) = \widehat{f}
\]
that vanishes in the past, aka the retarded solution.
Boundary conditions must not be given,
by causal independence forward in time.
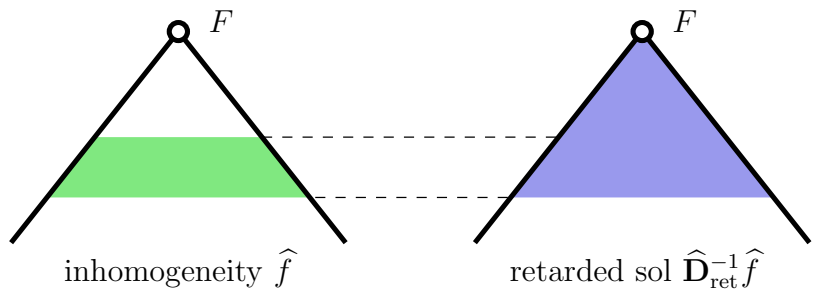
\begin{center}
\input{in.pstex_t}
\end{center}

A very simple scenario is \emph{finite codimension stability}:
That is, when there exist finitely many
functions $\widehat{u}_1,\ldots,\widehat{u}_J : \ncomp$
in the kernel of $\widehat{\op}$,
and complex valued linear functionals $\ell_1,\ldots,\ell_J$ such that,
for all nice $\widehat{f}$,
\[
\ret\widehat{f} \;-\;
\widehat{u}_1\,\ell_1(\widehat{f})
\;-\; \ldots \;-\;
\widehat{u}_J\,\ell_J(\widehat{f})
\]
decays
towards the future limit point $F$.
In this paper we identify
simple abstract assumptions on $\widehat{\op}$,
on top of the ones already stated,
that imply finite codimension stability.
\step
With a hyperbolic equation,
any non-smoothness in $\widehat{f}$
can persist in the retarded solution,
which can easily spoil finite codimension stability.
To avoid this, we put ourselves in function spaces
of {\infdiff} functions\footnote{%
These spaces are consistent with
nontrivial functions having compact support.
As for dropping $\infty$-differentiability,
consider Counterexample \ref{counterexample:nonsmooth}
and see Remark \ref{remark:nonsmooth}.
}.
\step
$\selfsim$-periodicity is the only symmetry assumption\footnote{%
In particular, spherical symmetry is not assumed.
Note that Choptuik's critical spacetime, a discretely self-similar
$1+3$ dim solution to Einstein's equations coupled to a
massless scalar field, is spherically symmetric.
An analogous solution to the vacuum Einstein equations,
if such exists, cannot be spherically symmetric.
}.
\step
This paper leaves open whether its results can be applied to
the linearized Einstein equations about
Choptuik's critical spacetime.
One could start from the proof of existence of this spacetime in
\cite{choptuikexists}.
Incidentally,
even though \cite{choptuikexists} is not about stability,
the proof of existence there
has some technical overlap with this paper,
see Example \ref{example:dim1}.

\subsection*{Intuition}

Let $\widehat{u}$ be a high frequency solution to
 the homogeneous equation $\widehat{\op}\widehat{u}=0$,
and pretend that the solution's wavelength
does not change much under propagation, as measured in the \intpic.
The same solution looks rather different in the \wrkpic,
where its wavelength increases:
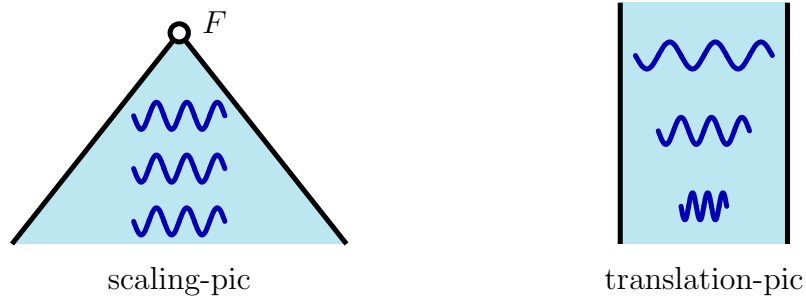
\begin{center}
\input{wave.pstex_t}
\end{center}

There is a more invariant, picture-independent way of saying this:
Pick a horizontal cross-section $B$, a ball,
take snapshots of the solution on
$\selfsim^p(B)$ for $p \in \nonnegint$, and pull them back to $B$.
These pullbacks
$\widehat{u}|_{\selfsim^p(B)} \circ \selfsim^p$
are a sequence of functions $B \to \C^N$
with increasing wavelengths.
\step
Clearly not all symmetric hyperbolic operators behave like this.
We impose a simple condition
on the principal part of $\widehat{\op}$
that gives such a regularizing effect at high frequencies;
this is the assumption later called \ref{s:def}.
\subsection*{Outline}

We use
coordinates
in which the cylinder is given by
\begin{align*}
\uc \;& =\; \mathrlap{\big\{ (\x{0},\ldots,\x{n}) \in \R^{1+n}
\;\big|\;
(\x{1})^2 + \ldots + (\x{n})^2 \leq 1
\big\}}\rule{90mm}{0pt}\\
\;&=\; \rule{0pt}{12pt} \mathrlap{\R \;\times\; \text{(closed unit ball in $\R^n$)}}
\end{align*}
and in which the self-diffeomorphism is given by
\[
\selfsim\;:\;
(\x{0},\x{1},\ldots,\x{n})
\;\mapsto\; (\x{0}+2\pi,\x{1},\ldots,\x{n})
\]
Many calculations in this paper are carried out
on the compact quotient
\begin{align*}
\dom \;& =\; \mathrlap{\uc \big/ \selfsim} \rule{90mm}{0pt}\\
\;& =\;
\rule{0pt}{12pt} \mathrlap{(\R/2\pi \allint) \;\times\; \text{(closed unit ball in $\R^n$)}}
\end{align*}

There are four abstract assumptions on the operator $\widehat{\op}$.
For the purpose of this introduction,
we state them informally:
\begin{enumerate}
\item[\ref{s:sh}]
$\selfsim$-periodicity and linear symmetric hyperbolicity,
with $\x{0}$ as time.
\item[\ref{s:bdy}]
Causal independence, forward in time.
\item[\ref{s:def}]
Condition on the principal part,
significant at high frequencies.
\item[\ref{s:ra}]
Regularity condition on the operator itself.
\end{enumerate}
\step

The main goal is to construct the right-inverse $\ret$,
aka the retarded Green's function,
and to show finite codimension stability.
By our earlier definition, we must show
the existence of
a finite-rank\footnote{%
Suppose $\banach,\banach'$ are vector spaces.
A linear operator $\banach \to \banach'$
is finite-rank iff
it is the composition of a linear map $\banach \to \C^J$ with
a linear map $\C^J \to \banach'$, for some integer $J$.
}
operator $\fr$ such that
$\widehat{\op}\circ \fr = 0$
and such that all elements in the image of
\[
\ret\;-\;\fr
\]
decay as $\x{0}\to+\infty$.
Finite codimension stability is \theoremref{fds}.
\step
\theoremref{fds} is proved by a contour integration argument.
This is based on
a Fourier-like transform,
used to write the retarded Green's function
as a contour integral in terms of the resolvent of
\[
\op
\;=\;
\widehat{\op}\big|_{\text{$\selfsim$-periodic functions}}
\]
where the vertical line means `restricted to'.
The two operators $\widehat{\op}$ and $\op$
are obviously equivalent, but it is always understood that:
\begin{align*}
\text{$\widehat{\op}$
maps functions $\uc \to \C^N$ to functions $\uc \to \C^N$}\\
\text{$\op$
maps functions $\dom \to \C^N$ to functions $\dom \to \C^N$}
\end{align*}
\step
Studying the resolvent of $\op$
is the main technical task.
This resolvent $(\op + z\mathbbm{1})^{-1}$
requires function spaces to make sense,
and may not exist for all $z\in \C$.
We study the resolvent on
suitable Banach spaces of
{\infdiff}
functions $\dom\to \C^N$
and show that
for $\RE z$ bigger than some constant $\rz \in \R$,
the resolvent exists and is compact\footnote{%
Suppose $\banach,\banach'$ are Banach spaces. A linear
map $\banach\to \banach'$
is compact iff every bounded sequence in $\banach$ is mapped
to a sequence in $\banach'$ that contains a Cauchy subsequence.
}.
By compactness at even just one point,
the resolvent
extends meromorphically to $\C$,
and the spectral projection associated to each pole
is finite-rank, that is, multiplicities are finite.
\step

The resolvent is periodic under $z \sim z + i$
up to conjugation\footnote{%
Just like the wave vector of a `Bloch wave'
on a lattice is defined only modulo
the reciprocal lattice.
Here the lattice is
$\x{0} \sim \x{0}+2\pi$,
the reciprocal lattice is $z\sim z+i$.
}.
In particular, $z$ is a pole if and only if $z+i$ is a pole.
A fundamental domain
such as $\{z \in \C \mid 0 \leq \IM z < 1\}$
decomposes into three pieces:
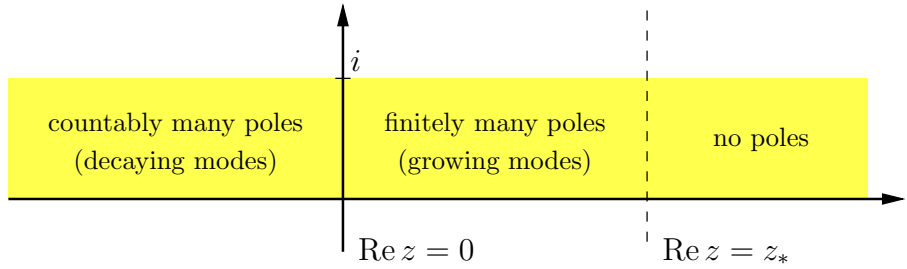
\begin{center}
\input{fin.pstex_t}
\end{center}

For every $c \in \R$, only finitely many poles
have $\RE z \geq c$ and  $0 \leq \IM z < 1$.
Therefore only finitely many
have $\RE z \geq 0$ and $0 \leq \IM z < 1$;
their spectral projections
yield an explicit formula for the finite-rank operator
$\fr$.
\step
One can think of individual poles
as decaying or growing modes,
but since the
spectral theorem for self-adjoint operators on Hilbert spaces
does not apply here,
we keep away from eigendecompositions.

\newcommand{\xs}{\x{1}_{\ast}}
\section{Examples and counterexamples}\label{sec:excex}

Before stating our abstract assumptions on $\op$,
we discuss examples.
They are invariant under arbitrary
translations in $\x{0}$, hence trivially
$\selfsim$-periodic.
\step
Let $\p_0,\ldots,\p_n$ be the partial
derivatives, $\p_ix^j = \delta_i^j$.
\begin{example}\label{example:dim1}
In $n=1$ and $N=1$ consider
\[
\op\;=\;\p_0 + \mu (\x{1}-\xs) \p_1
\]
with $\mu > 0$ and $-1 \leq \xs \leq 1$.
It satisfies all abstract assumptions.
Such operators are used in \cite{choptuikexists},
where the inverse of
$\p_0 + \mu (\x{1}+1)\p_1 + \mu$
is calculated explicitly,
on some space of real analytic functions
 on $\dom = (\R/2\pi \allint) \times [-1,1]$, using Fourier-Chebyshev series.
This inverse has regularizing features
at high frequencies,
used in \cite{choptuikexists}
to reduce the proof of existence of Choptuik's spacetime
to a finite, if very large, computer calculation for the low frequencies.
The regularizing features of this inverse
are similar to the estimates
on the resolvent of $\op$ obtained in this paper, in a more general setting.
\end{example}
\begin{example}\label{example:dim2}
In $n=3$ and $N=2$ consider
\[
\op\;=\;
\begin{pmatrix}
\phantom{\mu}\p_0 + \phantom{i}\mu \p_3 & \mu \p_1 + i\mu \p_2\\
\mu \p_1 - i\mu \p_2 & \phantom{\mu}\p_0 - \phantom{i}\mu \p_3
\end{pmatrix}
+ \begin{pmatrix} \mu & 0 \\ 0 & \mu \end{pmatrix}
\,(\x{1}\p_1 + \x{2}\p_2 + \x{3}\p_3)
\]
with $\mu > 0$.
It satisfies all abstract assumptions.
\end{example}

In Examples \ref{example:dim1} and \ref{example:dim2},
consider the eigenvalue problem
$\ker(\op + z\mathbbm{1}) \neq 0$ with $z\in \C$.
The point is that every eigenfunction $u: \dom \to \C^N$,
or rather its lift $\widehat{u}: \uc \to \C^N$,
yields a homogeneous solution,
 $\widehat{\op}(e^{z\x{0}}\widehat{u})=0$.
There are real analytic eigenfunctions of the form
\[
e^{iq\x{0}} P(\x{1},\ldots,\x{n})
\]
where $q \in \allint$, where
$P$ is an $N$-component polynomial of total degree $p \in \nonnegint$,
and the eigenvalue is $z = -iq - \mu p$.
In $P$
one can freely choose all terms of homogeneous degree $p$,
and then all terms of degree less than $p$ are recursively determined.
The set $i\allint - \mu \nonnegint$
of such eigenvalues is discrete, a half-lattice,
periodic under $z\sim z+i$,
and has real parts bounded from above.

\begin{counterexample}[eigenfunctions not \infdiff]\label{counterexample:nonsmooth}
The operator in
Example \ref{example:dim1}
has other eigenfunctions such as $\max\{0,\x{1}-\xs\}^s$
for all $s \in \C$ with say
$\RE s > 1$ and eigenvalue $z = -\mu s$.
The set of such eigenvalues is not discrete\footnote{%
Such eigenvalues can end up in the right half-plane,
say for operators of the form
$\op = \p_0 + \mu \x{1}\p_1 + \text{const}$,
with $\mu > 0$,
which also satisfy the abstract assumptions.}.
In this paper we avoid such eigenvalues by putting ourselves in
suitable Banach spaces of {\infdiff} functions.
\end{counterexample}

\begin{counterexample}[no causal independence]
Take Example \ref{example:dim1}
with $\xs < -1$.
It satisfies all abstract assumptions except
\ref{s:bdy}.
Since $\x{1} - \xs > 0$ everywhere on $\dom
= (\R/2\pi\allint)\times [-1,1]$,
there are real analytic eigenfunctions
$(\x{1} - \xs)^s$ for all $s \in \C$, with $z = -\mu s$.
\end{counterexample}

\begin{counterexample}[no high frequency effect]
Consider $\op = \p_0$. It satisfies all
abstract assumptions except \ref{s:def}.
It features infinite multiplicities
even for real analytic eigenfunctions.
\end{counterexample}

\section{Abstract assumptions}\label{sec:aa}

Let $n \geq 1$ be an integer.
Our working domain and its boundary are:
\begin{align*}
\dom
\;& =\;
\big\{
(\x{0},\x{1},\ldots,\x{n})
\in (\R/2\pi \allint) \times \R^n
\;\big|\;
(\x{1})^2 + \ldots + (\x{n})^2 \leq 1
\big\}\\
\p \dom
\;&=\;
\big\{
(\x{0},\x{1},\ldots,\x{n})
\in (\R/2\pi \allint) \times \R^n
\;\big|\;
(\x{1})^2 + \ldots + (\x{n})^2 = 1
\big\}
\end{align*}
Here $\x{0}\sim \x{0}+2\pi$.
The partial derivatives
$\p_0,\p_1,\ldots,\p_n$ satisfy
$\p_i x^j=\delta_i^j$.
Along $\p \dom$, the outward unit normal is
$(\omega_0,\omega_1,\ldots,\omega_n)
= (0,\x{1},\ldots,\x{n})$.
\step

\newcommand{\CN}[1]{\left|#1\right|_{\C^N}}
Given and fixed is a differential operator\footnote{%
As usual, $A^i\p_i = A^0\p_0 + \ldots + A^n\p_n$.
}\textsuperscript{,}\footnote{%
For simplicity, in Sections \ref{sec:intro} and \ref{sec:excex}
we used $\mathbbm{1}$
in places where we use $A^0$ from now on.
}
\[
\op\;=\; A^i\p_i + B
\]
subject to four abstract assumptions:
\begin{enumerate}[label=\textbf{(\roman*)}] 
\item\label{s:sh} $A^0,\ldots,A^n,B$ are $N \times N$ matrices
whose entries are {\infdiff} functions $\dom \to \C$.
The matrices $A^0,\ldots,A^n$ are Hermitian
and $A^0$ is positive definite ($\dagger$ is the
conjugate transpose):
\begin{align*}
(A^i)^{\dagger} \; & =\; A^i\\
A^0 \;& >\;0
\end{align*}
\item\label{s:bdy} $A^i\omega_i$ is positive semidefinite,
$A^i\omega_i \geq 0$,
along the boundary $\p\dom$.
\item\label{s:def} Set
\[
A^{ij} \;=\; \tfrac{1}{2}(\p_iA^j + \p_jA^i)
\]
There is a constant $\xi > 0$ and there are $N\times N$ matrices
$\auxi^0,\ldots,\auxi^n$ whose entries are {\infdiff} functions $\dom \to \C$ such that
\[
w_i^{\dagger}\big(
\xi A^{ij}
+
\tfrac{1}{2}A^i\auxi^j
+
\tfrac{1}{2}(\auxi^i)^{\dagger}A^j
\big) w_j
\;\geq\;
\delta^{ij}w_i^{\dagger}w_j
\]
for all $w_0,\ldots,w_n \in \C^N$ and
everywhere on $\dom$.
\item\label{s:ra} There is a sequence of constants
$\rr{0},\rr{1},\rr{2},\ldots \geq 0$
with
$\rr{0}=1$ and
\[
\rr{k+\ell} \;\leq\; \rr{k}\rr{\ell}
\]
for all $k,\ell \in \nonnegint$,
and there is a constant $Q > 0$ such that
\begin{align*}
\sum_{k=K+1}^{\infty}
(k+1)^{n/2}\;
\mathrlap{\frac{\rr{k-1} \na{k}}{k!}}
\rule{26mm}{0pt}
& \;\leq\; Q\rr{K}\\
\sum_{\;\;k=K\;\;}^{\infty}
(k+1)^{n/2}\;
\mathrlap{\frac{\rr{k} \nb{k}}{k!}}
\rule{26mm}{0pt}
& \;\leq\; Q\rr{K}(1+|z|)
\end{align*}
for $K=0,1$ (two values only) and for all $z\in \C$.
\end{enumerate}
In \ref{s:ra} we use the following supremum norms
of derivatives of order $k \in \nonnegint$:
\begin{align*}
\na{k}
\;&=\;
\sup_{\dom}
\sup_{\substack{\xx_{\alpha i} \in \C^N\\ \text{not all zero}}}
\frac{
\CN{
\sum_{|\alpha|=k}\sum_{i=0}^n \frac{k!}{\alpha!} (\p^{\alpha}A^i)\xx_{\alpha i}
}}{%
\big(
\sum_{|\alpha|=k}\sum_{i=0}^n
\frac{k!}{\alpha!} \CN{\xx_{\alpha i}}^2
\big)^{1/2}}\displaybreak[0]\\
\rule{0pt}{30pt}
\nb{k}
\;&=\;
\sup_{\dom}
\sup_{\substack{\xx_{\alpha} \in \C^N\\ \text{not all zero}}}
\frac{
\CN{%
\sum_{|\alpha|=k} \frac{k!}{\alpha!} (\p^{\alpha}(B+zA^0))\xx_{\alpha}
}}{%
\big(
\sum_{|\alpha|=k}
\frac{k!}{\alpha!} \CN{\xx_{\alpha}}^2
\big)^{1/2}}
\end{align*}
Throughout this paper, $\alpha \in \nonnegint^{1+n}$ 
is a multi-index, with
$|\alpha| = \alpha_0 + \ldots + \alpha_n$,
$\alpha! = \alpha_0!\cdots \alpha_n!$,
$\p^{\alpha} = (\p_0)^{\alpha_0}\cdots (\p_n)^{\alpha_n}$.
Also, $\CN{w} = (w^{\dagger}w)^{1/2}$ if $w \in \C^N$.
\step
Along with $\op$ itself, we consider
$\xi,\auxi^i,Q,\rr{\ell}$ to be part of the data.
That is, they are parameters of the theorems in this paper.
\begin{itemize}
\item
Some theorems require $\rr{\ell} > 0$ for all $\ell$.
\item
Some theorems impose a smallness condition on $\rr{1}$.
It is easy to satisfy both this smallness condition
and \ref{s:ra} simultaneously,
because of the following fact: If \ref{s:ra}
holds for a sequence $(\rr{\ell})$ then it continues to hold
for the sequence $(\kappa^\ell \rr{\ell})$ for all $0 < \kappa \leq 1$,
with the same $Q$.
\end{itemize}
\step

\emph{Informal discussion.}
Assumption \ref{s:sh} makes $\op$ 
symmetric hyperbolic,
in the sense of K.O.~Friedrichs.
Here $A^0>0$ holds uniformly, since $\dom$ is compact.
Assumption \ref{s:bdy}
is causal independence forward in time.
Assumption \ref{s:def} is a positivity condition for
the `deformation tensor' $A^{ij}$;
the naive condition with $\auxi^0=\ldots=\auxi^n=0$
is too strong, in fact
inconsistent
with \ref{s:sh}\footnote{%
Fix a nonzero $w \in \C^N$
and then set
$W = \int_0^{2\pi} \dd \x{0}\,w^{\dagger}A^{00}w$,
the integral taken along say $\x{1}=\ldots=\x{n} = 0$.
Assumption \ref{s:def} with $\auxi^0=\ldots=\auxi^n=0$
would imply $W > 0$,
whereas $A^{00} = \p_0A^0$ and periodicity
$\x{0}\sim \x{0} + 2\pi$ in \ref{s:sh} imply $W = 0$.
}.
\step

Assumption \ref{s:ra} are bounds for $A^i$ and $B$.
Two extreme cases are:
\begin{itemize}
\item $\rr{0}=1$ and $\rr{1} \geq 0$ but $\rr{2}=\rr{3}=\ldots=0$.
\item $\rr{\ell} = (\rr{1})^{\ell}$ for some $\rr{1} > 0$.
\end{itemize}
In the last case, $A^i$ and $B$ are real analytic.
The last case is extreme because
\ref{s:ra} requires
 $\rr{\ell} \leq (\rr{1})^{\ell}$.
Intermediate
cases include sequences with $\rr{\ell} > 0$
for all $\ell$ that go
to zero super-exponentially as $\ell \to \infty$.
\section{Main technical theorems, compact resolvent}\label{sec:maintechnical}

For all $u, v : \dom \to \C^N$
let $\ip{u}{v} = \int_{\dom} u^{\dagger} v$ be the standard inner product,
and let $\|u\| = \ip{u}{u}^{1/2}$.
Define the Sobolev seminorm of order $\ell \in \nonnegint$ by
\[
\|u\|_{\ell}\;=\;
\Big(\sum_{|\alpha|=\ell} \frac{\ell!}{\alpha!}\|\p^{\alpha}u\|^2\Big)^{1/2}
\]
We use the following function spaces\footnote{%
By definition, a function $\dom \to \C^N$ is {\infdiff}
if and only if it is the
restriction of an {\infdiff} function
$(\R/2\pi\allint)\times \R^n \to \C^N$.}\textsuperscript{,}\footnote{%
If $\rr{\ell}>0$ for all $\ell$, then $\ban{h}\subset \smoo$.
If $\rr{\ell}=0$ for one $\ell$ and therefore for almost all
$\ell$, then $\ban{h}$ is a Sobolev space
and $\smoo\subset\ban{h}$.
In either case,
$\smoo \cap \ban{h} \subset \ban{h}$ is dense.
}:
\begin{itemize}
\item
$\smoo
=
\{
u : \dom \to \C^N \mid
\text{$u$ is {\infdiff}}
\}$.
\item 
$\ban{h}= \{
u : \dom \to \C^N
\mid
\tn{u}_h < \infty
\}
$
for all $h \in \nonnegint$,
with norm
\[
\tn{u}_h\;=\;
\sum_{\ell=0}^{\infty}\frac{\rr{\ell}\|u\|_{\ell}}{(\ell+h)!}
\]
The Banach space $\ban{h}$
depends implicitly on the sequence $(\rr{\ell})$ in \ref{s:ra}.
\end{itemize}
We assume
\ref{s:sh}, \ref{s:bdy}, \ref{s:def}, \ref{s:ra}.
The theorems below refer to two constants\footnote{%
Their explicit values are in
Lemmas \ref{lemma:owhk} and \ref{lemma:wojxnks}.
}:
\begin{itemize}
\item A constant $\rz \in \R$ that depends only on $\op$.
\item A constant $\rrrz > 0$ that depends only on $\op,\xi,\auxi^i,Q$.
\end{itemize}
For each $z \in \C$ define\footnote{%
The operator $\op_z$ is used to construct
the retarded Green's function, in Section \ref{sec:fds}.
}
\[
\op_z \;=\; \op + z A^0
\]
\begin{theorem}\label{theorem:smoothbij}
If $\RE z \geq \rz$, then
$\op_z:\smoo\to\smoo$ is bijective.
\end{theorem}
\begin{theorem}\label{theorem:sjl}
If $\RE z \geq \rz$ and $\rr{1} \leq \rrrz$, then 
\[
\op_z^{-1}(\smoo\cap\ban{1})
\;\subset\;
\smoo\cap\ban{0}
\]
and it extends uniquely to a bounded
linear map $\op_z^{-1}: \ban{1} \to \ban{0}$.
\end{theorem}
\begin{theorem}\label{theorem:compactmeromorphic}
If $\RE z \geq \rz$
and $\rr{1} \leq \rrrz$
and $\rr{\ell}>0$ for all $\ell$, then
\[
\text{$\op_z^{-1}\;:\; \ban{1} \to \ban{1}$\;\; is compact}
\]
and there exists a unique extension to a meromorphic map
$w \mapsto \op_w^{-1}$ from $\C$ to the Banach space
of bounded linear maps $\ban{1}\to\ban{1}$.
This
extension satisfies the first resolvent identity
and is periodic up to conjugation\footnote{%
The three operators $A^0,\,e^{\pm i\x{0}}: \ban{1}\to\ban{1}$
are bounded
by Lemmas \ref{lemma:bstr}
and \ref{lemma:A0bounded}.}:
\begin{align*}
\op_{w}^{-1} - \op_{w'}^{-1}\;&=\; - (w-w') \op_w^{-1}A^0\op_{w'}^{-1}\\
\op_{w+i}^{-1}\;&=\;e^{-i\x{0}} \op_w^{-1} e^{i\x{0}}
\end{align*}
For all $w \in \C$ away from poles,
$\op_w^{-1}$ is compact.
\end{theorem}
\theoremref{smoothbij} follows from Corollary \ref{corollary:injective}
and \lemmaref{surjective}.
\theoremref{sjl} then follows from \lemmaref{wojxnks}.
\step
If $\rr{\ell} > 0$ for all $\ell$,
then the inclusion
 $\ban{0} \hookrightarrow \ban{1}$ is compact
by \lemmaref{cpctemb},
and hence \theoremref{sjl} implies compactness of $\op_z^{-1}$
in \theoremref{compactmeromorphic}.
In turn, compactness
implies the unique meromorphic extension:
it can be explicitly defined by fixing any $z \in \C$ with $\RE z \geq \rz$
and setting
\[
\op_w^{-1}
\;=\;
\op_z^{-1}\big(\mathbbm{1}+(w-z)A^0\op_z^{-1}\big)^{-1}
\]
which is meromorphic in $w \in \C$
by the spectral theory of compact operators,
applied to
 $A^0\op_z^{-1}: \ban{1}\to \ban{1}$.
This definition
is consistent when $\RE w \geq \rz$.

\section{Lemmas and proofs}
\begin{lemma}\label{lemma:epssub}
Suppose $\rr{\ell} > 0$ for all $\ell$.
Then for all $\eps > 0$,
every sequence in $\ban{0}$ with $\ban{0}$-diameter $\leq 1$
has a subsequence with $\ban{1}$-diameter $\leq \eps$.
Here the diameter is the supremum of all
pairwise distances.
\end{lemma}

\begin{PROOF}
Pick an integer $k \geq 1$ big enough
to make $\frac{1}{k+1} \leq \frac{\eps}{2}$.
By $\rr{\ell}>0$ for all $\ell$,
the given sequence in $\ban{0}$ is bounded in
the Sobolev space of order $k$.
Hence it has a Cauchy subsequence in the Sobolev space of order $k-1$,
by Rellich's theorem\footnote{%
The theorem applies because $\dom$ is compact.}.
In particular, it has a subsequence $(u_p)_{p\geq 0}$ such that
\begin{align*}
\sum_{\ell = 0}^{k-1}
 \frac{\rr{\ell}\|u_p-u_q\|_{\ell}}{(\ell+1)!}
 \;&\leq\;
 \frac{\eps}{2}
\intertext{for all $p,q$. On the other hand,}
\sum_{\ell = k}^{\infty}
 \frac{\rr{\ell}\|u_p-u_q\|_{\ell}}{(\ell+1)!}
 \;&\leq\;
 \frac{1}{k+1}\tn{u_p-u_q}_0
 \;\leq\;
 \frac{1}{k+1}
 \;\leq\;
 \frac{\eps}{2}
\end{align*}
for all $p,q$. Therefore $\tn{u_p-u_q}_1 \leq \eps$ for all $p,q$
as required.
\end{PROOF}
\begin{lemma}\label{lemma:cpctemb}
Suppose $\rr{\ell} > 0$ for all $\ell$.
Then $\ban{0}\hookrightarrow \ban{1}$ is compact.
\end{lemma}

\begin{PROOF}
Given
is a sequence in $\ban{0}$ with $\ban{0}$-diameter $\leq 1$.
A subsequence
with $\ban{1}$-diameter $\leq \eps$
will be called an $\eps$-subseq;
an $\eps$-subseq exists for all $\eps>0$ by \lemmaref{epssub}.
First pick a $\frac{1}{2}$-subseq.
Given this $\frac{1}{2}$-subseq,
store its first element and pick a $\frac{1}{3}$-subseq of the rest.
Given this $\frac{1}{3}$-subseq,
store its first element and pick a $\frac{1}{4}$-subseq of the rest.
And so forth.
The subsequence of stored elements
is a Cauchy sequence relative to $\ban{1}$.
\end{PROOF}
\begin{lemma}\label{lemma:owhk}
There exist constants $\rz \in \R$ and $\cc > 0$ such that
if $\RE z \geq \rz$ then for all
$u \in \smoo$ the functions
$J^i = \tfrac{1}{2} u^{\dagger}A^iu$ satisfy both
\begin{alignat*}{6}
\p_iJ^i & \;\leq\; &\; -J^0 \;& + \RE(u^{\dagger}\op_zu)\\
\p_iJ^i & \;\leq\; &\; -\cc_z (u^{\dagger}u) \;&+ \RE(u^{\dagger}\op_zu)
\end{alignat*}
everywhere on $\dom$, where
$\cc_z = \cc(1+|\RE z|)$.
Convention: All occurrences of $\rz$ and $\cc$
in this paper refer to their values as determined
by this lemma.
\end{lemma}

\begin{PROOF}
Using \ref{s:sh} we have
$\p_i J^i = -u^{\dagger}K_z u
 + \RE(u^{\dagger}\op_zu)$
where
\[
K_z \;=\; \tfrac{1}{2}(-\p_iA^i+B+B^{\dagger})\;+\;(\RE z)A^0
\]
is Hermitian.
Since $A^0 > 0$ uniformly,
we can choose $\rz$ such that $\RE z \geq \rz$
implies $K_z \geq \tfrac{1}{2}A^0$,
and choose $\cc>0$ such that $\RE z \geq \rz$ implies
$K_z \geq \cc_z \mathbbm{1}$.
\end{PROOF}

\begin{lemma}\label{lemma:sdhkkskka}
If $\RE z \geq \rz$ then
\[
\ip{u}{u}\;\leq\; \frac{1}{\cc_z} \RE\,\ip{u}{\op_z u}
\]
for all $u \in \smoo$.
\end{lemma}

\begin{PROOF}
By the divergence theorem and by \ref{s:bdy},
\[
\int_{\dom}\p_iJ^i \;=\; \int_{\p \dom}J^i\omega_i\;\geq\;0
\]
Now use the second inequality in \lemmaref{owhk}.
\end{PROOF}


\begin{corollary}\label{corollary:injective}
If $\RE z \geq \rz$ then
$\op_z: \smoo \to \smoo$ is injective.
\end{corollary}
The next lemma is where
$\p_jA^i$ comes out of hiding;
recall that its symmetric part is the deformation tensor
in \ref{s:def}.
The lemma splits the commutator
$[\p^{\alpha},\op_z]$
into two parts:
one part with derivatives of order $|\alpha|$,
the other part with derivatives of order less than $|\alpha|$,
aka lower order terms (lot).
\newcommand{\xxb}[2]{[\p^{#1},\op_{#2}]_{\textup{lot}}}
\newcommand{\xxbhat}[2]{[\p^{#1},\widehat{\op}_{#2}]_{\textup{lot}}}
\begin{lemma}\label{lemma:commutator}
Let $e_i \in \nonnegint^{1+n}$ be the $i$-th unit vector.
We have\footnote{%
It is clear from context where summation
over $i=0\ldots n$ and/or $j=0\ldots n$ is implicit.
}
\[
[\p^{\alpha},\op_z] \;=\;
\alpha_j (\p_jA^i) \p^{\alpha+e_i-e_j}
 + \xxb{\alpha}{z}
\]
where, by definition,
\[
\xxb{\alpha}{z}
 \;=\;
\sum_{\substack{\beta \leq \alpha \\|\beta| \geq 2}}
{\alpha \choose \beta}
(\p^{\beta}A^i)\p^{\alpha-\beta+e_i}
\;+\;
\sum_{\substack{\beta \leq \alpha \\|\beta| \geq 1}}
{\alpha \choose \beta}
(\p^{\beta}(B+zA^0))
\p^{\alpha-\beta}
\]
\end{lemma}

\begin{PROOF}
The product rule.
\end{PROOF}

\begin{lemma}\label{lemma:ags8sa}
For all $u \in \smoo$ and $\ell \in \nonnegint$ we have
\begin{align*}
& \RE \sum_{|\alpha|=\ell}\frac{\ell!}{\alpha!}
u^{\dagger}_{\alpha}([\p^{\alpha},\op_z]-\xxb{\alpha}{z})u\\
& \qquad \;\geq\;
\frac{\ell}{\xi} \sum_{|\alpha|=\ell}\frac{\ell!}{\alpha!}
u_{\alpha}^{\dagger}u_{\alpha}
\;-\; \frac{\ell}{\xi}
\RE \sum_{|\beta|=\ell-1}\frac{(\ell-1)!}{\beta!}
(\auxi^iu_{\beta+e_i})^{\dagger}(A^ju_{\beta+e_j})
\end{align*}
where $u_{\alpha} = \p^{\alpha}u$. This lemma uses notation from \ref{s:def}.
\end{lemma}

\begin{PROOF}
On the left hand side, use \lemmaref{commutator},
substitute $\alpha = \beta + e_j$,
use
$\alpha_j\frac{1}{\alpha!} = \frac{1}{\beta!}$
and
$\RE(u_{\beta+e_j}^{\dagger}(\p_jA^i)u_{\beta+e_i})
= u_{\beta+e_i}^{\dagger} A^{ij} u_{\beta+e_j}$;
the last because $A^i$ is Hermitian by
\ref{s:sh}.
One finds that the left hand side is equal to
\[
\ell\; \sum_{|\beta|=\ell-1}
\frac{(\ell-1)!}{\beta!}
u_{\beta+e_i}^{\dagger}A^{ij}u_{\beta+e_j}
\]
Now the lemma follows from \ref{s:def}
and the combinatorial \lemmaref{combinatorial}.
\end{PROOF}

\begin{lemma} \label{lemma:combinatorial}
Let $c_{\alpha} \in \C$ be a collection of
complex numbers, where $\alpha$ runs over $|\alpha|=\ell$
for some integer $\ell \geq 1$. Then
\[
\sum_{|\beta|=\ell-1}
\frac{(\ell-1)!}{\beta!}
\sum_{i=0}^n
c_{\beta+e_i}
\;=\;
\sum_{|\alpha|=\ell} \frac{\ell!}{\alpha!} c_{\alpha}
\]
\end{lemma}

\begin{PROOF}
It suffices to check this for
$c_{\alpha} = \delta_{\alpha = \gamma}$ for all fixed $|\gamma| = \ell$.
That is, all we need is
$\sum_{i=0}^n \frac{(\ell-1)!}{(\gamma-e_i)!}
= \frac{\ell!}{\gamma!}$,
which is true, and the lemma is proved.
\end{PROOF}

\begin{lemma}\label{lemma:cstr}
Let $a_{\alpha},b_{\alpha},c_{\alpha k} \in \smoo$
be collections of functions, where $\alpha$ runs over $|\alpha|=\ell
\in \nonnegint$, and where $k$ runs over a finite index set.
Then
\begin{align*}
\Big|\sum_{|\alpha|=\ell}
\frac{\ell!}{\alpha!}
\ip{a_{\alpha}}{b_{\alpha}}\Big|
\;&\leq\;
\Big(
\sum_{|\alpha|=\ell}\frac{\ell!}{\alpha!}
\|a_{\alpha}\|^2
\Big)^{1/2}
\Big(
\sum_{|\alpha|=\ell}\frac{\ell!}{\alpha!}
\|b_{\alpha}\|^2
\Big)^{1/2}\\
\Big(
\sum_{|\alpha|=\ell}\frac{\ell!}{\alpha!}
\|\sum_k c_{\alpha k}\|^2
\Big)^{1/2}
\;&\leq\;
\sum_k
\Big(
\sum_{|\alpha|=\ell}
\frac{\ell!}{\alpha!}
\|c_{\alpha k}\|^2
\Big)^{1/2}
\end{align*}
\end{lemma}
\begin{PROOF}
This is a Cauchy-Schwarz inequality and a triangle inequality.
\end{PROOF}

We now define a few abbreviations
that depend implicitly on a function
$u \in \smoo$ and on a number $z \in \C$.
As before, $u_{\alpha} = \p^{\alpha} u$.
For all $\ell \in \nonnegint$ set
\begin{align*}
\trd_{\ell} \;&=\;
\Big(
\sum_{|\alpha|=\ell}
\frac{\ell!}{\alpha!}
\|[\p^{\alpha},\op_z]u\|^2
\Big)^{1/2}\\
\fst_{\ell} \;&=\;
\Big(\sum_{|\alpha|=\ell}\frac{\ell!}{\alpha!}
\ltwo{\xxb{\alpha}{z}u}^2\Big)^{1/2}\\
\snd_{\ell} \;&=\;
\Big(\sum_{|\alpha|=\ell}\frac{\ell!}{\alpha!}
\ltwo{A^iu_{\alpha+e_i}}^2\Big)^{1/2}
\intertext{For all integers $0 \leq k \leq \ell$ set}
\mathcal{A}_{\ell k}
\;&=\; \Big(
\sum_{|\alpha|=\ell}\frac{\ell!}{\alpha!}
\Big\|
\sum_{\substack{\beta \leq \alpha\\ |\beta|=k}}
{\alpha \choose \beta}
(\p^{\beta}A^i)u_{\alpha-\beta+e_i}
\Big\|^2
\Big)^{1/2}\\
\mathcal{B}_{\ell k}
\;&=\; \Big(
\sum_{|\alpha|=\ell}\frac{\ell!}{\alpha!}
\Big\|
\sum_{\substack{\beta \leq \alpha\\ |\beta|=k}}
{\alpha \choose \beta}
(\p^{\beta}(B+zA^0))u_{\alpha-\beta}
\Big\|^2
\Big)^{1/2}
\end{align*}
\begin{lemma} \label{lemma:ell}
If $\RE z \geq \rz$, then for all $u \in \smoo$ and $\ell \in \nonnegint$
we have
\[
\|u\|_{\ell}
\;\;\leq\;\;
\frac{\xi}{\xi \cc_z + \ell}
\big(
\|\op_z u\|_{\ell}
+ \fst_{\ell}
\big)
\;+\;
\frac{\ell |\auxi|_0}{\xi \cc_z + \ell} \snd_{\ell-1}
\]
Here $|\auxi|_0$
is defined just like $\na{0}$ is defined in Section \ref{sec:aa}.
\end{lemma}


\begin{PROOF}
Set $f = \op_z u$
and
$f_{\alpha} = \p^{\alpha}f$. Trivially,
$\op_z u_{\alpha}=f_{\alpha} - [\p^{\alpha},\op_z]u$.
Since $\RE z \geq \rz$, \lemmaref{sdhkkskka} implies
\[
\cc_z \ip{u_{\alpha}}{u_{\alpha}}\;\leq \;
\RE \ip{u_{\alpha}}{f_{\alpha} - [\p^{\alpha},\op_z]u}
\]
Equivalently,
\[
\cc_z \ip{u_{\alpha}}{u_{\alpha}}
+ \RE \ip{u_{\alpha}}{
([\p^{\alpha},\op_z]-\xxb{\alpha}{z})
u}
\;\leq \;
\RE \ip{u_{\alpha}}{f_{\alpha} - \xxb{\alpha}{z}u}
\]
Multiply both sides by $\ell!/\alpha!$
and then sum over $|\alpha| = \ell$.
On the left hand side, use \lemmaref{ags8sa}.
Then use \lemmaref{cstr} in a number of places to get
\[
\Big(\cc_z + \frac{\ell}{\xi}\Big) \|u\|_{\ell}^2
\;\leq\;
\|u\|_{\ell}\,\big(\|f\|_{\ell} + \fst_{\ell}\big)
+
\frac{\ell}{\xi}
|\auxi|_0 \|u\|_{\ell} \snd_{\ell-1}
\]
where, for the rightmost term, one also needs
\[
\Big(
\sum_{|\beta|=\ell-1}\frac{(\ell-1)!}{\beta!}
\|\auxi^iu_{\beta+e_i}\|^2
\Big)^{1/2}
\;\leq\;
|\auxi|_0\|u\|_{\ell}
\]
which follows from the definition of $|\auxi|_0$
and
\lemmaref{combinatorial}.
\end{PROOF}

\begin{lemma} \label{lemma:triangle}
We have $\snd_{\ell}\;\leq\;
\|\op_zu\|_{\ell}
+
\trd_{\ell}
+
\mathcal{B}_{\ell 0}$.
\end{lemma}

\begin{PROOF}
Use the triangle inequality in \lemmaref{cstr},
together with
\begin{align*}
A^iu_{\alpha+e_i}
& = (\op_z-B-zA^0)u_{\alpha}\\
& = \p^{\alpha}(\op_z u) - [\p^{\alpha},\op_z]u
- (B+zA^0)u_{\alpha}
\end{align*}
\end{PROOF}

\begin{lemma}\label{lemma:triangletwo}
We have\footnote{%
The summation
$\sum_{k=K}^{\ell}$ gives zero when $\ell < K$.
}
\begin{align*}
\trd_{\ell} + \mathcal{B}_{\ell 0} & \;\leq\;
\sum_{k=1}^{\ell} \mathcal{A}_{\ell k}
+
\sum_{k=0}^{\ell} \mathcal{B}_{\ell k}\\
\fst_{\ell} & \;\leq\;
\sum_{k=2}^{\ell} \mathcal{A}_{\ell k}
+
\sum_{k=1}^{\ell} \mathcal{B}_{\ell k}
\end{align*}
\end{lemma}

\begin{PROOF}
Use the triangle inequality in \lemmaref{cstr}.
More in detail, every summation $\sum_{\beta \leq \alpha}$
as in \lemmaref{commutator}
is written as $\sum_k \sum_{\beta \leq \alpha, |\beta|=k}$,
and then the triangle inequality is used for the summation over $k$.
\end{PROOF}


\begin{lemma} \label{lemma:summary}
If $\RE z \geq \rz$, then for all $u \in \smoo$ and $\ell \in \nonnegint$
we have
\begin{align*}
\|u\|_{\ell} \;& \leq\;
\frac{\xi}{\xi \cc_z +\ell}
 \Big(
\|\op_z u\|_{\ell}
+ 
\sum_{k=2}^{\ell}
\mathcal{A}_{\ell k}
+ 
\sum_{k=1}^{\ell}
\mathcal{B}_{\ell k}
\Big)\\
& \qquad + \frac{\ell|\auxi|_0}{\xi \cc_z + \ell}
 \Big(
\|\op_z u\|_{\ell-1}
+
\sum_{k=1}^{\ell-1}
\mathcal{A}_{\ell-1, k}
+ 
\sum_{k=0}^{\ell-1}
\mathcal{B}_{\ell-1, k}
\Big)
\end{align*}
\end{lemma}
\begin{PROOF}
Lemmas \ref{lemma:ell}, \ref{lemma:triangle}, \ref{lemma:triangletwo}.
The special case $\ell = 0$, which
simplifies to $\|u\| \leq \cc_z^{-1} \|\op_z u\|$,
is also a direct corollary of \lemmaref{sdhkkskka}.
\end{PROOF}

\begin{lemma} \label{lemma:shkhfdls}
For all integers $0 \leq k \leq \ell$,
\begin{align*}
\mathcal{A}_{\ell k} & \;\leq\;
{\ell \choose k}^{1/2} {\ell + n \choose k}^{1/2}
\,\na{k}\, \|u\|_{\ell-k+1}\\
\mathcal{B}_{\ell k} & \;\leq\;
{\ell \choose k}^{1/2} {\ell + n \choose k}^{1/2}
\,\nb{k}\, \|u\|_{\ell-k}
\end{align*}
\end{lemma}

\begin{PROOF}
The definitions of $\mathcal{A}_{\ell k}$ and $\na{k}$ imply
\[
\mathcal{A}_{\ell k}
\;\leq\;
\na{k}\; \Big(
\sum_{|\alpha|=\ell}\frac{\ell!}{\alpha!}
\sum_{\substack{\beta \leq \alpha\\|\beta|=k}}
\sum_{i=0}^n
\frac{\beta!}{k!}
{\alpha \choose \beta}^2
\| u_{\alpha-\beta+e_i}\|^2
\Big)^{1/2}
\]
It is convenient to replace
$\|u_{\alpha-\beta+e_i}\|^2 = \sum_{|\gamma| = \ell-k+1} \|u_{\gamma}\|^2 \delta_{\alpha-\beta+e_i = \gamma}$ and to move the summation over $\gamma$
to the left, which gives
\[
\mathcal{A}_{\ell k}
\;\leq\;
\na{k}\;\Big(\sum_{|\gamma|=\ell-k+1}
c_{\gamma}
\frac{(\ell-k+1)!}{\gamma!} \|u_{\gamma}\|^2
\Big)^{1/2}
\]
with the purely combinatorial coefficient
\[
c_{\gamma}\;=\;
\frac{\gamma!}{(\ell-k+1)!}
\sum_{|\alpha|=\ell}\frac{\ell!}{\alpha!}
\sum_{\substack{\beta \leq \alpha\\|\beta|=k}}
\sum_{i=0}^n
\frac{\beta!}{k!}
{\alpha \choose \beta}^2
\delta_{\alpha-\beta+e_i = \gamma}
\]
One can check that
$c_{\gamma} = {\ell \choose k} {\ell + n \choose k}$
if $|\gamma|=\ell-k+1$, which is independent of $\gamma$, and
this gives the desired estimate
for $\mathcal{A}_{\ell k}$.
Similar for $\mathcal{B}_{\ell k}$.
\end{PROOF}

\begin{lemma} \label{lemma:conv}
For $K=0,1$ we have
\begin{align*}
\sum_{\ell = 0}^{\infty}
\frac{\rr{\ell}}{(\ell+1)!}
\sum_{k=K+1}^{\ell}
\mathcal{A}_{\ell k}
\;&\leq\;
Q\rr{K} \tn{u}_0\\
\sum_{\ell = 0}^{\infty}
\frac{\rr{\ell}}{\ell!}
\sum_{k=K}^{\ell}
\mathcal{B}_{\ell k}
\;&\leq\;
Q\rr{K}(1+|z|) \tn{u}_0
\end{align*}
\end{lemma}

\begin{PROOF}\rule{0pt}{0pt}
\begin{itemize}
\item
For the 1st estimate, use
$\rr{\ell} \leq \rr{k-1} \rr{\ell-k+1}$ where $1\leq k \leq \ell$.\\
For the 2nd estimate, use $\rr{\ell} \leq \rr{k} \rr{\ell-k}$
where $0\leq k \leq \ell$.
\item Use \lemmaref{shkhfdls}
with
${\ell \choose k}^{1/2}
{\ell+n \choose k}^{1/2}
\leq {\ell \choose k} (k+1)^{n/2}$
 when $0 \leq k \leq \ell$.
\item
For the 1st estimate, use
$\frac{1}{(\ell+1)!} {\ell \choose k} \leq \frac{1}{k! (\ell-k+1)!}$.\\
For the 2nd estimate, use
$\frac{1}{\ell!} {\ell \choose k} \leq \frac{1}{k!(\ell-k)!}$.
\item
$\sum_{\ell=0}^{\infty}
\sum_{k=K}^{\ell} c_{k,\ell-k}
= \sum_{k=K}^{\infty}\sum_{p=0}^{\infty}c_{k,p}$
for all $c: \nonnegint \times \nonnegint \to [0,\infty)$.
\item 
Use assumption \ref{s:ra}.
\end{itemize}
\end{PROOF}

\begin{lemma}\label{lemma:wojxnks}
Define $\rrrz > 0$ by
\[
Q\rrrz \big(
\xi + 3 \cc^{-1}
 + |\auxi|_0 + 2\cc^{-1}\xi^{-1}|\auxi|_0
\big)\;=\;\tfrac{1}{2}
\]
If $\RE z \geq \rz$ and $\rr{1} \leq \rrrz$ then,
for all $u \in \smoo$,
\[
\tn{u}_0 \;\leq\;
2 e^{2\rr{1} |\IM z|}
\big(
\xi + \cc^{-1}
+ |\auxi|_0 \rr{1} \big)\,
\tn{\op_z u}_1
\]
\end{lemma}

\begin{PROOF}
We show that the inequality holds without the exponential factor
when $|\IM z| \leq 1$;
the general case then follows from
$\op_{z+i}=e^{-i\x{0}} \op_z e^{i\x{0}}$ and \lemmaref{bstr}.
\lemmaref{summary} implies
\begin{align*}
\|u\|_{\ell} \;& \leq\;
(\xi + \cc^{-1})
\frac{1}{\ell+1}
 \Big(
\|\op_z u\|_{\ell}
+ 
\sum_{k=2}^{\ell}
\mathcal{A}_{\ell k}\Big)
+
\frac{1}{\cc_z}
\sum_{k=1}^{\ell}
\mathcal{B}_{\ell k}
\\
& \qquad + |\auxi|_0
 \Big(
\|\op_z u\|_{\ell-1}
+
\sum_{k=1}^{\ell-1}
\mathcal{A}_{\ell-1, k}\Big)
+ 
\frac{\ell|\auxi|_0}{\xi \cc_z}
\sum_{k=0}^{\ell-1}
\mathcal{B}_{\ell-1, k}
\end{align*}
Multiply both sides by $\rr{\ell}/\ell!$ and then sum over $\ell$.
\lemmaref{conv}
and the fact that
$(1+|z|)\cc_z^{-1} \leq 2 \cc^{-1}$,
because we are assuming $|\IM z|\leq 1$, imply
\begin{align*}
\tn{u}_0 &\;\leq\;
(
\xi + \cc^{-1}
)
\big(
\tn{\op_z u}_1
+
Q\rr{1} \tn{u}_0
\big)
+ 2\cc^{-1} Q\rr{1} \tn{u}_0\\
&\qquad +
|\auxi|_0
\big(
\rr{1} \tn{\op_z u}_1
+
Q\rr{1} \tn{u}_0
\big)
+
2\cc^{-1} \xi^{-1} |\auxi|_0 Q\rr{1} \tn{u}_0
\end{align*}
using $\rr{\ell} \leq \rr{1} \rr{\ell-1}$
and $\rr{0}=1$.
The assumption $\rr{1} \leq \rrrz$ makes
the coefficient of $\tn{u}_0$
on the right hand side $\leq \frac{1}{2}$.
Solve for $\tn{u}_0$ and be done\footnote{%
Actually,
$\tn{u}_0$ could be infinite,
and then the proof does not work.
However,
if $\rr{\ell}=0$ for one $\ell$ and hence almost all $\ell$,
then $\tn{u}_0<\infty$
and the proof does work.
The general case then follows from
truncating the sequence $(\rr{\ell})$
at some index and taking the limit.}.
\end{PROOF}
\begin{lemma} \label{lemma:bstr}
The two operators $e^{\pm i\x{0}}: \ban{h} \to \ban{h}$
are bounded for all $h \in \nonnegint$.
Explicitly,
$\tn{e^{\pm i\x{0}} u}_{h}
\leq
e^{\rr{1}} \tn{u}_{h}$
for all $u \in \smoo$.
\end{lemma}

\begin{PROOF}
By routine estimation, $\|e^{\pm i\x{0}}u\|_{\ell}
\leq
\sum_{m=0}^{\ell} {\ell \choose m}
\|u\|_{\ell-m}$.
Then use the definition of $\tn{\,\cdot\,}_h$ and
$\frac{1}{(\ell+h)!}{\ell \choose m}
\leq \frac{1}{m!(\ell-m+h)!}$ and
$\rr{\ell} \leq (\rr{1})^m \rr{\ell-m}$.
\end{PROOF}

\begin{lemma}\label{lemma:A0bounded}
The operator $A^0: \ban{h}\to \ban{h}$ is bounded for all $h \in \nonnegint$.
\end{lemma}

\begin{PROOF}
Note that
$\|(B+zA^0)u\|_{\ell}\leq
\sum_{k=0}^{\ell}
\mathcal{B}_{\ell k}$.
Similar to \lemmaref{conv},
\[
\sum_{\ell = 0}^{\infty}
\frac{\rr{\ell}}{(\ell+h)!}
\sum_{k=0}^{\ell}
\mathcal{B}_{\ell k}
\;\leq\;
Q(1+|z|) \tn{u}_h
\]
using $\frac{1}{(\ell+h)!}{\ell \choose k}
\leq \frac{1}{k!(\ell-k+h)!}$.
Hence
$\tn{(B+zA^0)u}_h \leq Q(1+|z|)\tn{u}_h$.
Hence the map
$B+zA^0: \ban{h}\to\ban{h}$ is bounded for all $z\in \C$,
hence $A^0$ is.
\end{PROOF}


\begin{lemma}\label{lemma:decay}
Let $\uc_{\x{0} \geq 0} \subset \uc$ be the subset
of the universal cover space corresponding to $\x{0}\geq 0$.
Suppose
$\RE z \geq \rz$ and suppose
\[
\widehat{u}\;:\; \uc_{\x{0} \geq 0} \;\to\; \C^N
\]
is {\infdiff} and satisfies $\widehat{\op}_z \widehat{u} = 0$.
Then each derivative $\widehat{u}_{\alpha}
= \p^{\alpha}\widehat{u}$
converges to zero exponentially fast as $\x{0} \to +\infty$,
uniformly in $\x{1},\ldots,\x{n}$.
\end{lemma}

\begin{PROOF}
Set
\[
E_{\ell}(\x{0})
\;=\;
\sum_{|\alpha|=\ell}\frac{\ell!}{\alpha!}
\int_{\text{unit ball in $\R^n$}}
\big(\tfrac{1}{2}
\widehat{u}^{\dagger}_{\alpha} A^0 \widehat{u}_{\alpha}
\big)_{\text{at time $\x{0}$}}
\]
It suffices to show that
for all $\ell \in \nonnegint$
and all $0 \leq c < 1$,
\[
\textstyle\sup_{\x{0}\geq 0}
|e^{c\x{0}} E_{\ell}(\x{0})| \;<\; \infty
\]
The proof is by induction over $\ell$.
Since $\RE z \geq \rz$,
the first inequality in \lemmaref{owhk}\footnote{%
This lemma on $\dom$ also holds on the universal cover $\uc$.}
and the divergence theorem and \ref{s:bdy} imply
\[
\frac{\dd}{\dd \x{0}}
E_{\ell}(\x{0})
\;\leq\;
- E_{\ell}(\x{0})
+ \RE 
\sum_{|\alpha|=\ell}\frac{\ell!}{\alpha!}
\int_{\text{unit ball in $\R^n$}}
\big(
\widehat{u}_{\alpha}^{\dagger}
\widehat{\op}_z \widehat{u}_{\alpha}
\big)_{\text{at time $\x{0}$}}
\]
In the second term on the right hand side,
replace\footnote{%
Recall from \ref{s:sh} that $A^i,B$ are \infdiff,
hence the commutator is defined.
}
\[
\widehat{\op}_z\widehat{u}_{\alpha}
\; = \;
-([\p^{\alpha},\widehat{\op}_z]
  -\xxbhat{\alpha}{z})\widehat{u}
-\xxbhat{\alpha}{z}\widehat{u}
\]
Now use \lemmaref{ags8sa}\footnote{%
This lemma on $\dom$ also holds on the universal cover $\uc$.
Recall that it relies on \ref{s:def}.
}.
The first term coming from
\lemmaref{ags8sa}
has favorable sign and is dropped.
In the second term coming from
\lemmaref{ags8sa}
we replace $A^j\widehat{u}_{\beta+e_j}
= [\widehat{\op}_z,\p^{\beta}]\widehat{u}
- (B+zA^0) \widehat{u}_{\beta}$ where always $|\beta|=\ell-1$.
Now triangle and Cauchy-Schwarz inequalities imply
\[
\tfrac{\dd}{\dd \x{0}}
E_{\ell}
\;\leq\;
- E_{\ell}
+ \text{(const)}\;
(E_0^{1/2} + \ldots + E_{\ell-1}^{1/2}) E_{\ell}^{1/2}
\]
given that $A^0>0$ holds uniformly by \ref{s:sh}.
The unspecified constant depends on many things,
including $\ell$,
but it does not depend on $x \in \uc$.
We get
\[
\tfrac{\dd}{\dd \x{0}}
(e^{c\x{0}/2} E_{\ell}^{1/2})
\;\leq\; \text{(const)} e^{c\x{0}/2}
(E_0^{1/2} + \ldots + E_{\ell-1}^{1/2})
\]
whenever $0\leq c < 1$.
By the induction hypothesis,
the integral of the right hand side over $\x{0} \in [0,\infty)$
is finite, and the lemma follows.
\end{PROOF}

\begin{lemma}\label{lemma:surjective}
If $\RE z \geq \rz$
then $\op_z: \smoo \to \smoo$ is surjective.
\end{lemma}

\begin{PROOF}
Given any $f \in \smoo$, we construct a
 $u \in \smoo$ such that
$\op_z u = f$.
Let $\widehat{f}: \uc \to \C^N$
be the lift of $f$ to the universal cover.
Let
\[ \widehat{u}\;:\; \uc_{\x{0} \geq 0} \;\to\; \C^N \]
be the unique {\infdiff} solution to
$\widehat{\op}_z\widehat{u} = \widehat{f}$
that vanishes at $\x{0} = 0$\footnote{%
Existence and uniqueness for linear symmetric
hyperbolic systems.
Use \ref{s:sh}
and \ref{s:bdy}.}.
Since $\widehat{f}\circ \selfsim = \widehat{f}$
we have
$\widehat{\op}_z(\widehat{u} - \widehat{u}\circ \selfsim) = 0$.
By $\RE z \geq \rz$ and
\lemmaref{decay} we have
$\widehat{u} - \widehat{u} \circ \selfsim \to 0$
exponentially fast
as $\x{0} \to +\infty$;
same for all derivatives of all orders.
Hence
$\widehat{u}\circ \selfsim^p - \widehat{u} \circ \selfsim^{p+1} \to 0$
exponentially fast as $p \to \infty$
on the compact
$\uc_{0 \leq \x{0} \leq 4\pi}$.
This being exponentially fast,
we have
$\widehat{u}\circ \selfsim^p - \widehat{u} \circ \selfsim^q \to 0$ as $p,q \to \infty$.
The limit of $\widehat{u}\circ \selfsim^p$
on $\uc_{0 \leq \x{0} \leq 4\pi}$
descends to an {\infdiff}
$u: \dom \to \C^N$. Since $\widehat{\op}_z(\widehat{u}\circ \selfsim^p)
= \widehat{f}$ for all $p$, we have
$\op_z u = f$.
\end{PROOF}

\newcommand{\fin}{p}
\newcommand{\supernice}{\banc{-\infty}{\infty}}
\section{Finite codimension stability}\label{sec:fds}
This section relies heavily on \theoremref{compactmeromorphic}.
To be able to freely use this theorem, we always assume
in this section:
\begin{center}
$\rr{1} \leq \rrrz$ and $\rr{\ell} > 0$ for all $\ell$
\end{center}
This assumption is repeated explicitly in \theoremref{fds} below,
for emphasis,
but the assumption is implicit throughout this section.
\step
Define the following function spaces:
\begin{itemize}
\item
$\smooc
=
\{
\widehat{u} : \uc \to \C^N \mid
\text{$\widehat{u}$ is {\infdiff}}
\}$.
\item
For all fixed $a < b$, denote by $\banc{a}{b}$ the
set of all $\widehat{u} \in \smooc$ such that
\[
\sum_{\ell=0}^{\infty}\frac{\rr{\ell}}{(\ell+1)!}
\Big(
\sum_{|\alpha|=\ell} \frac{\ell!}{\alpha!}
\int_{\uc}
e^{-2c\x{0}} |\p^{\alpha}\widehat{u}|^2
\Big)^{1/2}\;<\;\infty
\]
for each $c \in \R$ with $a < c < b$.
\end{itemize}

Roughly, a function is in $\supernice$ if it is
smooth as given by the sequence $(\rr{\ell})$,
and if it and all its derivatives of all orders 
decay super-exponentially for both $\x{0}\to \pm \infty$.
This space allows us to state a clean theorem;
one can certainly relax the assumptions of the theorem in various directions.

\begin{theorem}[Finite codimension stability]\label{theorem:fds}
In addition to the abstract assumptions in Section \ref{sec:aa},
assume $\rr{1} \leq \rrrz$ and $\rr{\ell} > 0$ for all $\ell$.
Let
\[
\mathllap{\ret}\;:\;\supernice\;\to\;\smooc
\]
be the retarded Green's function of $\widehat{\op}$,
explicitly constructed in \lemmaref{retgreensfunction}.
There exists a finite-rank operator\footnote{%
If $\banach$ and $\banach'$ are vector spaces,
then a linear operator $\banach \to \banach'$
is a finite-rank operator iff
it is the composition of a linear map $\banach \to \C^J$ with
a linear map $\C^J \to \banach'$, for some $J < \infty$.
The maps should be continuous,
but for simplicity, we do not introduce topologies.
}
\[
\mathllap{\fr}\;:\;\supernice\;\to\;\smooc
\]
such that $\widehat{\op}\circ \fr = 0$
and such that all elements in the image of
\[
\ret\;-\;\fr
\]
decay exponentially fast as $\x{0}\to+\infty$,
uniformly in $\x{1},\ldots,\x{n}$,
and the same decay statement holds for all
partial derivatives of all orders.
\end{theorem}
In the remainder of this section, we prove this theorem
and provide details, including an explicit
formula for $\fr$
and a description of its image.
\begin{lemma}\label{lemma:trivialdecay}
For all $\widehat{u} \in \banc{a}{b}$:
\begin{itemize}
\setlength{\itemsep}{1pt}
\item If $a < 0$ then $\widehat{u}$ decays exponentially fast as $\x{0} \to +\infty$.
\item
If $b > 0$ then $\widehat{u}$ decays exponentially fast as $\x{0} \to -\infty$.
\item
If $a =-\infty$ then $\widehat{u}$ decays super-exponentially fast as $\x{0} \to +\infty$.
\item
If $b = +\infty$ then $\widehat{u}$ decays super-exponentially fast as $\x{0} \to -\infty$.
\end{itemize}
The same decay statements hold for all partial derivatives of all orders
of $\widehat{u}$.
\end{lemma}

Define $S_{a,b} = \{z \in \C \mid a < \RE z < b \}$,
an infinite vertical strip.
\begin{lemma}[Fourier-like transform] \label{lemma:flt}
In the following, let $u_z$ be the evaluation of the map
$u_{\ast}$ at the point $z \in S_{a,b}$.
There is a linear, bijective map
\[
\left\{
u_{\ast} : S_{a,b} \to \ban{1}
\;\left|\;
\begin{aligned}
&\text{$u_{\ast}$ is holomorphic}\\
&\text{$u_{z+i} = e^{-i\x{0}} u_z$ for all $z \in S_{a,b}$}
\end{aligned}
\right.\right\}
\;\to\,
\banc{a}{b}
\]
given by $u_{\ast} \mapsto \widehat{u}$ where
\[
\widehat{u}(\x{}) \;=\;
\frac{1}{i}\int_{z'}^{z'+i}\dd z\; e^{z\x{0}} u_z(\x{})
\]
for all $x \in \uc$, all $z'\in S_{a,b}$ and all paths contained in $S_{a,b}$.
Convention: From here on, it is implicit that $u_{\ast}$
denotes the transform
of $\widehat{u}$ and conversely.
\end{lemma}
\begin{PROOF}
The integral does not
depend on the choice of path,
because
the integrand is holomorphic in $z \in S_{a,b}$ and periodic under 
$z \sim z+i$.
For injectivity, use a straight path $t \mapsto c + it$
with $a < c < b$ to get
\[
\widehat{u}(\x{0}+2\pi \fin,\x{1},\ldots,\x{n})
\;=\; \textstyle\int_0^1 \dd t\,
 e^{(c+it)\x{0}} e^{(c+it) 2\pi \fin } u_{c+it}(x)
\]
for all $\fin \in \allint$,
and now note that if $\widehat{u}=0$
then $\R/\allint \to \C,
\,t \mapsto e^{(c+it)\x{0}}u_{c+it}(x)$
is a function all whose Fourier coefficients vanish,
hence $u_{\ast}=0$.
For surjectivity, one checks that 
every $\widehat{u} \in \banc{a}{b}$
is the image point of $u_{\ast}$ given by
\[
u_z \;=\; \textstyle\sum_{\fin \in \allint}
(\,e^{-z\x{0}} \widehat{u}\,) \circ \selfsim^p
\]
Apart from estimates that we omit\footnote{%
Useful lemma: For all $z \in \C$ and $\widehat{u} \in \smooc$
and all `weight' functions $w: \uc \to [0,\infty)$:
\begin{multline*}
\sum_{\ell=0}^{\infty}
\frac{\rr{\ell}}{(\ell+1)!}
\Big(
\sum_{|\alpha|=\ell} \frac{\ell!}{\alpha!}
\int_{\uc} w\,|\p^{\alpha}(e^{z\x{0}}\widehat{u})|^2
\Big)^{1/2}\\
\;\leq\;
e^{\rr{1}|z|}
\sum_{\ell=0}^{\infty}
\frac{\rr{\ell}}{(\ell+1)!}
\Big(
\sum_{|\alpha|=\ell} \frac{\ell!}{\alpha!}
\int_{\uc} w\,e^{2(\RE z)\x{0}} |\p^{\alpha}\widehat{u}|^2
\Big)^{1/2}
\end{multline*}\par}, this concludes the proof.
\end{PROOF}

We now define $z_{\ast\ast\ast} \leq
z_{\ast\ast} \leq \rz$ with $z_{\ast\ast\ast} < 0$ by
\begin{align*}
z_{\ast\ast}\;&=\;
\sup\,\{ \RE z \mid \text{$z$ is a pole}\}\\
z_{\ast\ast\ast}\;&=\;
\sup\,\{ \RE z \mid \text{$z$ is a pole with $\RE z < 0$}\}
\end{align*}
where the poles are those of the map
$z\mapsto \op_z^{-1}$ in \theoremref{compactmeromorphic}\footnote{%
We write $\sup$ (rather than $\max$)
only because there could be no poles at all.}.
The set of poles is periodic under $z \sim z+i$.
Here is an example:
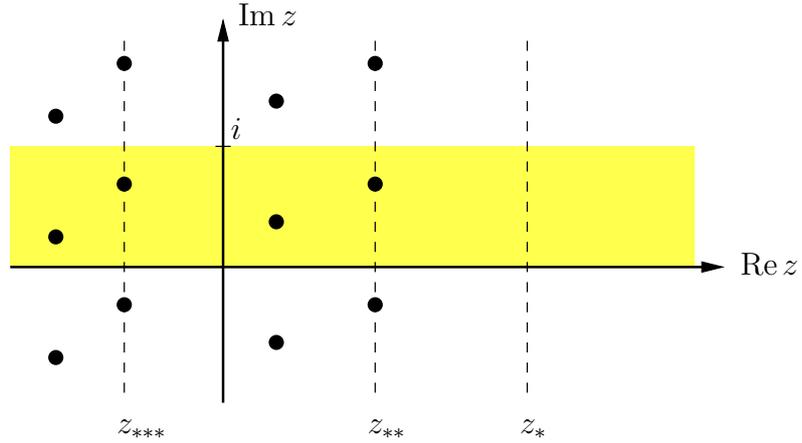
\begin{center}
\input{ctr0.pstex_t}
\end{center}
The domain indicated in this figure,
$\{z \in \C \mid 0 \leq \IM z < 1\}$, is a fundamental domain for $z\sim z+i$.
Here exactly two poles have nonnegative real
part in a fundamental domain,
but in general there could be any finite number.
\begin{lemma}\label{lemma:abc}
For all $z\in \C$ away from the poles of $z \mapsto \op_z^{-1}$,
and for all $f \in \ban{1}$,
applying $\widehat{\op}$
to the function
$x\mapsto e^{z\x{0}} (\op_z^{-1} f)(x)$ gives 
$x\mapsto e^{z\x{0}} f(x)$.
\end{lemma}

\begin{PROOF}
If $\RE z \geq \rz$ then $\op_z^{-1}$
is an honest inverse
in the sense of \theoremref{smoothbij}
and the claim follows from the operator identity
$\widehat{\op}e^{z\x{0}} = e^{z\x{0}}(\widehat{\op} + zA^0)$.
This implies the claim for general $z$ by a meromorphic
continuation argument.
\end{PROOF}

\begin{lemma}\label{lemma:retgreensfunction}
Define a linear map
\[
\ret\;:\; \supernice \;\to\;\banc{z_{\ast\ast}}{\infty}
\]
by
$\widehat{f}\mapsto \widehat{u}$
where $u_z = \op^{-1}_z f_z$, using \lemmaref{flt}. Explicitly,
\[
(\ret \widehat{f})(x)\;=\;
\frac{1}{i}\int_{z'}^{z'+i}
\dd z\, e^{z\x{0}}\, (\op^{-1}_zf_z)(x)
\]
for all $x \in \uc$
and all paths contained in $S_{z_{\ast\ast},\infty}$.
Then this is a right-inverse of $\widehat{\op}$,
and more specifically, it is the retarded Green's function.
\end{lemma}
\begin{PROOF}
It is a right-inverse by \lemmaref{abc}.
The map $\ret$ is the retarded Green's function
because all elements in its image converge to zero
super-exponentially as $\x{0}\to -\infty$,
and so do all partial derivatives of all orders.
\end{PROOF}

We now define
\[
\Lambda\;=\;\{z \mid \text{$z$ is a pole with $\RE z \geq 0$}\}
\;\cap\; \{z \mid 0 \leq \IM z < 1\}
\]
The second factor is a fundamental domain for $z\sim z+i$,
and any other fundamental domain could be used instead.
This is always a finite set, \mbox{$|\Lambda| < \infty$}.
\begin{lemma}\label{lemma:xxx1}
Define
$\fr:\supernice\to\smooc$
by
\[
(\fr \widehat{f})(x)
\;=\;
\sum_{\lambda \in \Lambda}
\frac{1}{i}\int_{\textup{loop about $\lambda$}}
\dd z\, e^{z\x{0}}\,(\op_z^{-1}f_z)(x)
\]
for all $x \in \uc$. Then
\begin{align*}
\widehat{\op}\circ \fr \;&=\;0\\
\image(\ret-\fr)\;&\subset\;\banc{z_{\ast\ast\ast}}{0}
\end{align*}
Note: The decay in \theoremref{fds} now follows
from $z_{\ast\ast\ast}<0$ and \lemmaref{trivialdecay}.
\end{lemma}

\begin{PROOF}
For $\widehat{\op}\circ \fr = 0$
use \lemmaref{abc} and
$\int_{\textup{loop about $\lambda$}} \dd z\,e^{z\x{0}}
f_z(x) = 0$.
For the image, the case $z_{\ast\ast} < 0$ is trivial, because
on the one hand $z_{\ast\ast} = z_{\ast\ast\ast}$
and hence $\image(\ret) \subset \banc{z_{\ast\ast\ast}}{0}$,
and on the other hand $|\Lambda|=0$ and $\fr = 0$.
The case $z_{\ast\ast} \geq 0$ is conveniently
discussed using our $|\Lambda|=2$ example:
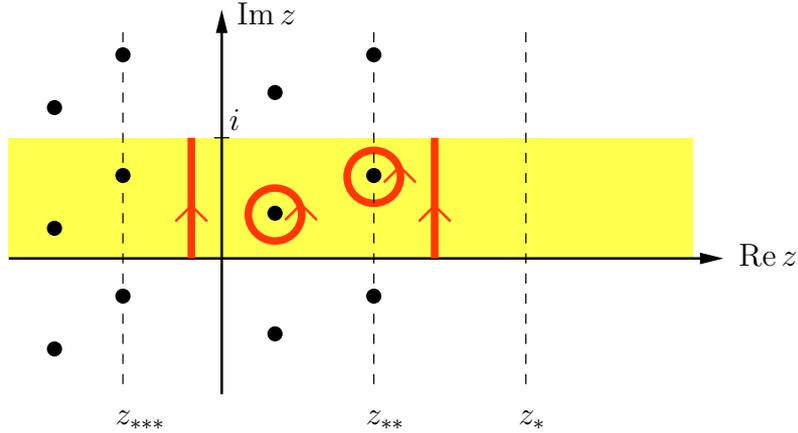
\begin{center}
\input{ctr1.pstex_t}
\end{center}
Let $\widehat{f} \in \supernice$. We want to show that
$(\ret - \fr)\widehat{f} \in \banc{z_{\ast\ast\ast}}{0}$.
 Fix $x\in \uc$.
Then $I_z(x) = i^{-1}
e^{z\x{0}}\,(\op_z^{-1}f_z)(x)$ is meromorphic in $z$
with poles coming from $\op_z^{-1}$ only,
and periodic under $z \sim z+i$. Observe that:
\begin{itemize}
\item
The integral of $I_z(x)$ along the path on the right gives
$(\ret \widehat{f})(x)$.
\item
The integral of $I_z(x)$ about the $|\Lambda|$-many loops gives
$(\fr \widehat{f})(x)$.
\end{itemize}
Their difference is equal to the integral of $I_z(x)$
along the path on the left, by Cauchy's theorem,
which yields an element of
$\banc{z_{\ast\ast\ast}}{0}$
by \lemmaref{flt}.
\end{PROOF}


\begin{lemma}\label{lemma:xx2}
If $\lambda \in \Lambda$
and $\ell \in \nonnegint$ then
\begin{align*}
(P_{\lambda \ell}f)(x)
\;&=\; \frac{1}{2\pi i}
\int_{\textup{loop about $\lambda$}}
\dd z\,(z-\lambda)^\ell\,(\op_z^{-1}f)(x)
\intertext{is a finite-rank operator
$P_{\lambda \ell}:\ban{1}\to\ban{1}$, and}
(\widehat{P}_{\lambda \ell}f)(x)
\;&=\;
\frac{1}{2\pi i}
\int_{\textup{loop about $\lambda$}}
\dd z\,(z-\lambda)^\ell\,e^{z\x{0}}\,(\op_z^{-1}f)(x)
\end{align*}
is a finite-rank operator
$\widehat{P}_{\lambda \ell} : \ban{1} \to \smooc$. We have
\begin{align*}
\widehat{\op}
\circ \widehat{P}_{\lambda \ell}\;&=\; 0\\
\image(\widehat{P}_{\lambda \ell}) \;&\subset\;
\textup{(polynomials in $\x{0}$)}\,
e^{\lambda \x{0}}\, \ban{1}
\end{align*}
\end{lemma}

\begin{PROOF}
We use \theoremref{compactmeromorphic}.
Since $\op_z^{-1}$ is compact away from poles,
$P_{\lambda \ell}$ is compact\footnote{%
The compact operators are a closed subspace of the space of bounded operators.
}. By the first resolvent identity,
$P_{\lambda k}A^0P_{\lambda \ell} = P_{\lambda(k+\ell)}$.
It follows that
$\image(P_{\lambda \ell}) \subset \image(P_{\lambda 0}A^0)$.
It also follows that $P_{\lambda 0}A^0$ 
is a projection.
As a compact projection,
$\dim(\image(P_{\lambda 0}A^0)) < \infty$.
Hence
$\dim(\image(P_{\lambda \ell})) < \infty$,
as claimed.
Now use
$e^{z \x{0}} =
\sum_{k=0}^{\infty}
\frac{1}{k!}  (z-\lambda)^k (\x{0})^k e^{\lambda \x{0}}$
which implies
\[
\widehat{P}_{\lambda \ell}
\;=\;
\sum_{k=0}^{\infty} \frac{1}{k!} (\x{0})^k e^{\lambda \x{0}}\,
P_{\lambda(k+\ell)}
\]
The sum is finite, because $P_{\lambda k} = 0$
if $k$ is equal to or bigger than the order of the pole at $\lambda$,
and therefore $\widehat{P}_{\lambda \ell}$ is a finite-rank operator
as well.
\end{PROOF}

\begin{lemma}
For all $\widehat{f} \in \supernice$ we have
\[
\fr \widehat{f}\;=\;
2\pi \sum_{\lambda \in \Lambda}
\sum_{\ell=0}^{\infty}
\frac{1}{\ell!}\, \widehat{P}_{\lambda \ell}
\,\Big(\frac{\dd^\ell f_z}{\dd z^{\ell}}\Big)_{\textup{at $z=\lambda$}}
\]
There are only finitely many pairs
$(\lambda,\ell) \in \Lambda \times \nonnegint$
for which $\widehat{P}_{\lambda \ell} \neq 0$.
\end{lemma}
\begin{PROOF}
This follows from Lemmas \ref{lemma:xxx1} and \ref{lemma:xx2},
using a Taylor expansion of
the holomorphic function $f_{\ast}$ about
each $\lambda \in \Lambda$.
\end{PROOF}

\begin{corollary}
The operator $\fr: \supernice \to \smooc$ is finite-rank and
\begin{align*}
\widehat{\op}\circ \fr \;&=\; 0\\
\image(\fr) \;&\subset\; \textstyle\sum_{\lambda \in \Lambda}
\textup{(polynomials in $\x{0}$)} e^{\lambda \x{0}} \ban{1}
\end{align*}
\end{corollary}
\section{Remarks}

\begin{remark}
The assumption that we are on the closed unit ball in
the $n$ spatial coordinates is not essential.
One can use
\[
\dom
\; =\;
\big\{
(\x{0},\x{1},\ldots,\x{n})
\in (\R/2\pi \allint) \times \R^n
\;\big|\;
(\x{1},\ldots,\x{n}) \in \Gamma(\x{0})
\big\}
\]
where $\Gamma(\x{0})\subset \R^n$ satisfies
$\Gamma(\x{0}+2\pi) = \Gamma(\x{0})$, all sufficiently smooth.
\end{remark}

\begin{remark}
This concerns another generalization that would require a lot of
careful checking.
Introduce an auxiliary positive definite (Riemannian) metric
$h = h_{ij}\dd x^i \otimes \dd x^j$ on $\dom$.
This paper may correspond to
the special case
$h_{ij} = \delta_{ij}$ of a more general,
more geometric formulation of the assumptions and theorems
that is invariant under diffeomorphisms of $\dom$,
with the understanding that $h$ is also transformed.
More general, because one now has the additional freedom of choosing $h$.
One can also try to introduce an auxiliary inner product on $\C^N$ etc.
None of this has been tried.
\end{remark}


\begin{remark}
In applications to general relativity,
one has to deal with gauge freedom.
One can gauge-fix and make the equations symmetric hyperbolic,
but it would also be interesting to try to reformulate the assumptions
and theorems in this paper in a more gauge-invariant way.
\end{remark}

\newcommand{\lz}{L}
\begin{remark} \label{remark:nonsmooth}
Despite
Counterexample \ref{counterexample:nonsmooth},
one can generalize
\theoremref{compactmeromorphic}
to finite differentiability.
One may not get a compact resolvent,
and it may not be meromorphic on $w \in \C$,
but it will be on half-planes $\RE w > w_{\ast}$;
the question is how negative one can take $w_{\ast}$
as a function of $\lz = \sup\{\ell \mid \rr{\ell}>0\}$.
The following sketch relies heavily on the operator seminorm
$\semi{\,\cdot\,}$ discussed separately in Appendix \ref{app:seminorm}.
Lemmas
\ref{lemma:epssub} and
\ref{lemma:cpctemb} generalize to
\[
\semi{\ban{0} \hookrightarrow \ban{1}} \;=\; \frac{1}{\lz+1}
\]
which is proved using Rellich's theorem. Informally: the bigger $\lz$,
the closer the inclusion operator is to being compact.
Using \lemmaref{wojxnks},
the first part of
 \theoremref{compactmeromorphic} generalizes to:
If $\RE z \geq \rz$ and $\rr{1} \leq \rrrz$
then
\[
\semi{\op_z^{-1}: \ban{1}\to \ban{1}}
\;\leq\;
\frac{
2 e^{2\rr{1} |\IM z|}
(\xi + \cc^{-1}
+ |\auxi|_0 \rr{1})
}{\lz+1}
\]
The numerator is written out for clarity;
here we consider the situation where all parameters 
in the numerator are fixed.
For $\lz=\infty$
we recover the original \theoremref{compactmeromorphic},
but we are now interested in $\lz < \infty$.
Putting things together, including Appendix \ref{app:seminorm},
one finds that
the bigger $\lz$, the more negative one can take $w_{\ast}$.
As $\lz \to \infty$
one can take $w_{\ast}\to -\infty$.
\end{remark}

\section*{Acknowledgments}
Thanks to Eugene Trubowitz for his interest and many
useful suggestions.

\appendix
\section{%
Non-compactness seminorm}\label{app:seminorm}

For a bounded operator $T$ on a Banach space, the operator norm
tells one that the resolvent
$(T-\lambda)^{-1}$
exists on $|\lambda| > \|T\|$.
By contrast, the operator seminorm defined below
tells one that the resolvent
exists on $|\lambda| > \semi{T}$, except for
a discrete set of points
that are actual eigenvalues.
\step
The spectral theory of compact operators is obtained as a special case,
because $\semi{T}=0$ iff $T$ is compact.
In fact, the arguments in this appendix are minor
adaptations
of standard arguments in
Riesz's spectral theory of compact operators.
Everything in this appendix
is probably
available in the literature on
`measures of non-compactness'.
\step
For all Banach spaces $\banach,\banach'$
and all bounded linear $T: \banach \to \banach'$, set
\[
\setlength{\abovedisplayskip}{15pt}
\setlength{\belowdisplayskip}{15pt}
\semi{T}\,=\,
\inf \left\{
\eps > 0
\left|
\begin{aligned}
&\text{\footnotesize For every sequence $(v_p)_{p\geq 0}$ in $\banach$
with diameter $\leq 1$,}\\
&\text{\footnotesize there is a subsequence of
$(Tv_p)_{p \geq 0}$ in $\banach'$ with diameter $\leq \eps$.}
\end{aligned}
\right.
\right\}
\]
Here the diameter is the supremum of all
pairwise distances. Then:
\begin{itemize}
\setlength{\itemsep}{2pt}
\item $0 \leq \semi{T} \leq \|T\|$.
\item $\semi{\,\cdot\,}$
is a continuous seminorm on the space of bounded operators.
\item $\semi{T}=0$ if and only if $T$ is compact.
\item $\semi{T_1 T_2} \leq \semi{T_1}\semi{T_2}$.
\item
If $\dim \banach=\infty$ then
$\semi{\mathbbm{1}_\banach} = 1$.
\end{itemize}
\step

For the remainder of this appendix,
fix a Banach space $\banach$ and a bounded linear
operator $T:\banach \to \banach$.
Define $\pspek \subset \spek \subset \C$ by:
\begin{align*}
\pspek
\;&=\;
\{\lambda \in \C \mid \text{$T-\lambda$ is not injective}\}\\
\spek\;&=\;
\{ \lambda \in \C \mid \text{$T-\lambda$ is not bijective}\}
\end{align*}
For all $\lambda \notin \spek$, the resolvent $(T-\lambda)^{-1}$
is a bounded operator, by the open mapping theorem.
The spectrum $\spek$ is a compact subset of $\C$.
\begin{lemma}\label{lemma:skhksad}
If $|\lambda| > \semi{T}$
then $\image(T-\lambda)$ is a closed subspace.
\end{lemma}
\begin{PROOF}
First observe that $\image(T-\lambda) = \image (T-\lambda)|_C$ where
\[
C\;=\;
\{
\,
v \in \banach \mid \dist (v, \ker(T-\lambda))
\geq \tfrac{1}{2} \|v\|
\,
\}
\]
Hence it suffices to show that,
if $(v_p)$ is a sequence in $C$ such that
$((T-\lambda)v_p)$ converges,
then $(v_p)$ has a Cauchy subsequence.
We combine\footnote{%
The second by
$\|(Tv_q-Tv_p) - \lambda(v_q-v_p)\| \to 0$
and the reverse triangle inequality.
}:
\begin{itemize}
\setlength{\itemsep}{2pt}
\item The definition of $\semi{T}$.
\item 
$\|Tv_q-Tv_p\| - |\lambda| \|v_q-v_p\|\,\to\, 0$
as $q,p\to\infty$.
\end{itemize}
They imply that \emph{there exists
a $0 < \kappa < 1$ such that
every subsequence of $(v_p)$ with diameter $d$
contains a subsequence with diameter $\leq \kappa d$}.
In fact, every $\kappa$ with $\semi{T}/|\lambda| < \kappa < 1$ will do.
This implies that
$(v_p)$ has a Cauchy subsequence,
and that we are done,
if $(v_p)$ has even just one subsequence that has finite diameter
($\Leftrightarrow$ that is bounded).
The remaining case is
$\|v_p\| \to \infty$.
Then $((T-\lambda)v_p/\|v_p\|)$ converges to zero.
Hence $(v_p/\|v_p\|)$ has a Cauchy subsequence,
by the argument just given,
with limit in $\ker(T-\lambda)$.
But this contradicts $v_p \in C$,
namely
$\dist (v_p/\|v_p\|,\ker(T-\lambda)) \geq \tfrac{1}{2}$,
hence $\|v_p\| \not\to \infty$.
\end{PROOF}

\begin{lemma} \label{lemma:subsp}
If $r > \semi{T}$ then there does not exist a
sequence $(\lambda_p)_{p\geq 0}$ of complex numbers,
and a sequence $(\banach_p)_{p\geq 0}$ of subspaces of $\banach$, such that:
\begin{itemize}
\setlength{\itemsep}{2pt}
\item $|\lambda_p| \geq r$.
\item $\banach_p$ is a closed subspace.
\item
$\banach_{p\pm 1} \subset \banach_p$ with proper inclusion,
and $(T-\lambda_p)\banach_p \subset \banach_{p\pm 1}$.
\end{itemize}
The sign $\pm$ is arbitrary,
but it is understood to be the same in both places.
\end{lemma}
\begin{PROOF} We prove the minus version;
the plus version is similar.
Suppose, by contradiction, that such sequences do exist.
By Riesz's lemma, there exist unit vectors $v_p \in \banach_p$ with
$\dist(\banach_{p-1},v_p) \geq \frac{1}{2}$.
Also note that $T\banach_p \subset \banach_p$.
For all integers $0 \leq q<p$ and $m\geq 1$ we have
\[
\|T^mv_q - T^mv_p\|
\;=\;
\|
\underbrace{T^mv_q}_{\in \banach_{p-1}}
- \underbrace{(T^m-\lambda_p^m)v_p}_{\in \banach_{p-1}}
  - \lambda_p^mv_p\|
\;\geq\;
\tfrac{1}{2}|\lambda_p|^m
\;\geq\;
\tfrac{1}{2}r^m
\]
The sequence $(v_p)$ has diameter $\leq 2$,
whereas every subsequence of $(T^mv_p)$
has diameter $\geq \tfrac{1}{2} r^m$.
Therefore $\semi{T^m} \geq \tfrac{1}{4}r^m$.
Therefore $\semi{T} \geq (\tfrac{1}{4})^{1/m} r$.
Since this holds for all $m$, we get $\semi{T} \geq r$,
a contradiction.
\end{PROOF}

\begin{lemma}
If $|\lambda| > \semi{T}$
and $\lambda \in \spek$ then $\lambda \in \pspek$.
\end{lemma}
\begin{PROOF}
The subspace $\banach_p = \image (T-\lambda)^p$ is closed,
by induction,
using $\banach_{p+1} = (T-\lambda)\banach_p$
and $\semi{T|_{\banach_p}} \leq \semi{T} < |\lambda|$
and \lemmaref{skhksad}.
Suppose, by contradiction, that
$\lambda \in \spek \setminus \pspek$.
Then $\banach_{p+1}\subset \banach_p$ is proper.
This contradicts the plus version of \lemmaref{subsp}
with $r = |\lambda|$
and $\lambda_p = \lambda$.
\end{PROOF}


\begin{lemma}
If $r > \semi{T}$
then $\{\lambda \in \pspek \mid |\lambda| \geq r\}$
is a finite set.
\end{lemma}

\begin{PROOF}
Suppose, by contradiction, that there exists a
sequence $(\lambda_p)_{p \geq 0}$ of pairwise distinct
$\lambda_p \in \pspek$ with $|\lambda_p| \geq r$.
Pick eigenvectors $v_p\neq 0$ with $Tv_p = \lambda_pv_p$
and set $\banach_p = \SPAN\{v_0,\ldots,v_p\}$. Then
$\dim \banach_p = p+1$,
since the $\lambda_p$ are distinct.
This contradicts the minus version of \lemmaref{subsp}.
\end{PROOF}


\newcommand{\httpref}[1]{\href{http://#1}{#1}}


{\footnotesize
\begin{thebibliography}{}
\bibitem{choptuik1993} Choptuik M.W.,
\href{http://link.aps.org/doi/10.1103/PhysRevLett.70.9}{%
Phys.~Rev.~Lett., 70, 9-12 (1993)}\\
\emph{Universality and scaling in gravitational collapse
of a massless scalar field}
\bibitem{GM} Gundlach C.~and Martin-Garcia J.M.,
\href{http://www.livingreviews.org/lrr-2007-5}{%
Living Rev.~Relativity 10 (2007), 5}\\
\emph{Critical Phenomena in Gravitational Collapse}
\bibitem{choptuikexists}
Reiterer M.~and Trubowitz E.,
\httpref{arxiv.org/abs/1203.3766}\\
\emph{Choptuik's critical spacetime exists}
\bibitem{gradedliealgebra}
Reiterer M.~and Trubowitz E.,
\httpref{arxiv.org/abs/1412.5561}\\
\emph{The graded Lie algebra of general relativity}
\end{thebibliography}}
\end{document}